# On *E*-Orlicz Theory

Abdulhameed Qahtan Abbood Altai and Nada Mohammed Abbas Alsafar

*Abstract*—In this paper, based on concepts of *E*-convex sets and *E*-convex functions, we introduce new concepts of *E*-*N*-functions, *E*-Young functions, *E*-strong Young functions and *E*-Orlicz functions by relaxing definitions of *N*-functions, Young functions, strong Young functions, Orlicz functions respectively. We also give concepts of *E*-Orlicz spaces, weak *E*-Orlicz spaces, *E*-Orlicz-Sobolev spaces, weak *E*-Orlicz-Sobolev spaces, *E*-Orlicz-Morrey spaces and weak *E*-Orlicz-Morrey spaces, *E*-Orlicz-Lorentz spaces and weak *E*-Orlicz-Lorentz spaces. In addition, we consider the implicit properties based on the effect of an operator *E* on the functions and the spaces.

## I. INTRODUCTION

BIRNBAUM and Orlicz introduced the Orlicz spaces in 1931 as a generalization of the classical Lebesgue spaces, where the function $u^p$ is replaced by a more general convex function $\Phi$ [2]. The concept of $E$-convex sets and $E$-convex functions were introduced by Youness to generalize the classical concepts of convex sets and convex functions to extend the studying of the optimality for non-linear programming problems in 1999 [3]. Chen defined the semi-$E$-convex functions and studied its basic properties in 2002 [3]. The concepts of pseudo $E$-convex functions and $E$-quasiconvex functions and strictly $E$-quasiconvex functions were introduced by Syau and Lee in 2004 [6]. The concept of Semi strongly $E$-convex functions was introduced by Youness and Tarek Emam in 2005 [8]. Sheiba Grace and Thangavelu considered the algebraic properties of $E$-convex sets in 2009 [4]. $E$-differentiable convex functions was defined by Meghed, Gomma, Youness and El-Banna [5] to transform a non-differentiable function to a differentiable function in 2013. Semi-$E$-convex function was introduced by Ayache and Khaled in 2015 [1].

Abdulhameed Qahtan Abbood Altai is with Department of Mathematics and Computer Science, University of Babylon, College of Basic Education, Babil, 51002, Iraq (ahbabil1983@gmail.com).

Nada Mohammed Abbas Alsafar is with the Department of Mathematics, College of Education for Pure Sciences, University of Babylon, Babil, 51002, Iraq (nadaalsafar333@gmail.com).

In our paper, our goal is to define the $E$-$N$-functions, $E$-Young functions, $E$-strong Young functions and $E$-Orlicz functions using the concepts of $E$-convex sets and $E$-convex functions to generalize and extend studying of the Orlicz theory via defining $E$-Orlicz spaces, weak $E$-Orlicz spaces, $E$-Orlicz-Sobolev spaces, weak $E$-Orlicz-Sobolev spaces, $E$-Orlicz-Morrey space, weak $E$-Orlicz-Morrey space, $E$-Orlicz-Lorentz spaces and weak $E$-Orlicz-Lorentz spaces generated by non-Young functions but $E$-Young functions with a map $E$.

**Contents of the paper**. The definitions of $E$-$N$-function, $E$-Young function, $E$-strong Young function and $E$-Orlicz function are presented in section II. The elementary properties of $E$-$N$-functions, $E$-Young functions, $E$-strong Young functions and $E$-Orlicz functions and their relationships will be considered in section III and IV respectively. And $E$-Orlicz spaces, weak $E$-Orlicz spaces, $E$-Orlicz-Sobolev spaces, weak $E$-Orlicz-Sobolev spaces, $E$-Orlicz-Morrey spaces, weak $E$-Orlicz-Morrey spaces, $E$-Orlicz-Lorentz spaces and weak $E$-Orlicz-Lorentz spaces will be stated in section V. In addition, the implicit properties will be established.

## II. PRELIMINARIES

**Definition 1[6].** A set $S \subset R^n$ is said to be $E$-convex iff there is a map $E: R^n \to R^n$ such that $\lambda E(x) + (1-\lambda)E(y) \in S$, for each $x, y \in S, 0 \leq \lambda \leq 1$.

**Definition 2[6].** A function $f: R^n \to R$ is said to be $E$-convex on a set $S \subset R^n$ iff there is a map $E: R^n \to R^n$ such that $S$ is an $E$-convex set and
$$f(\lambda E(x) + (1-\lambda)E(y)) \leq \lambda f(E(x)) + (1-\lambda)f(E(y)),$$
for each $x, y \in M$ and $0 \leq \lambda \leq 1$. And $f$ is called $E$-concave on a set $S$ if
$$f(\lambda E(x) + (1-\lambda)E(y)) \geq \lambda f(E(x)) + (1-\lambda)f(E(y)),$$
for each $x, y \in S$ and $0 \leq \lambda \leq 1$.

**Definition 3.** Let $S$ be a nonempty subset of $R^n$. A function $f: S \to R^m$ is said to be $E$-continuous at $a \in S$ iff there is a map $E: R^n \to R^n$ such that for every $\varepsilon > 0$ there is $\delta > 0$ implies
$$\|f(E(x)) - f(E(a))\| < \varepsilon \text{ whenever } \|x - a\| < \delta$$
and $f$ is said to be $E$-continuous on $S$ iff $f$ is $E$-continuous at every $x \in S$.

**Definition 4.** Let $(\Omega, \Sigma, \mu)$ be a measure space. A function $\Phi: \Omega \times [0, \infty) \to \mathbb{R}$ is called an $E$-$N$-function if there exists a map $E: \Omega \times [0, \infty) \to \Omega \times [0, \infty)$ such that for $\mu$-a.e. $t \in \Omega$, $[0, \infty)$ is an $E$-convex and $\Phi$ is an $E$-even, $E$-continuous, $E$-convex of $u$ on $[0, \infty)$, $\Phi(E(t,u)) > 0$ for any $u \in (0, \infty)$,
$$\lim_{u \to 0^+} \frac{\Phi(E(t,u))}{u} = 0, \lim_{u \to \infty} \frac{\Phi(E(t,u))}{u} = \infty$$
and for each $u \in [0, \infty)$, $\Phi(E(t,u))$ is an $\mu$-measurable function of $t$ on $\Omega$.

**Remark 5.** Every $N$-function is an $E$-$N$-function if the map $E$ is taken as the identity map. But it is not every $E$-$N$-function an $N$-function.

**Example 6.** Let $\Phi: \mathbb{R} \times [0, \infty) \to \mathbb{R}$ be defined by $\Phi(t,u) = tu^2$ and let $E: \mathbb{R} \times [0, \infty) \to \mathbb{R} \times [0, \infty)$ be defined as $E(t,u) = (|t|, u)$. Then $\Phi$ is an $E$-$N$-function but it is not an $N$-function because that for $\mu$-a.e. $t \in \mathbb{R}$, $\Phi(t,u)$ is concave of $u$ for $t \in (-\infty, 0)$.

**Example 7.** Let $\Phi: \mathbb{R} \times [0, \infty) \to \mathbb{R}$ be defined by $\Phi(t,u) = (1-t)u^2 + t \exp(u)$ and let $E: \mathbb{R} \times [0, \infty) \to \mathbb{R} \times [0, \infty)$ be defined by $E(t,u) = (t, \ln u^2)$. Then, $\Phi$ is an $E$-$N$-function but it is not an $N$-function since for $\mu$-a.e. $t \in \mathbb{R}$, $\Phi(t,u)$ is not even.

**Definition 8.** Let $(\Omega, \Sigma, \mu)$ be a measure space. A function $\Phi: \Omega \times [0, \infty) \to \mathbb{R}$ is called an $E$-Young function if there exists a map $E: \Omega \times [0, \infty) \to \Omega \times [0, \infty]$ such that for $\mu$-a.e. $t \in \Omega$, $[0, \infty)$ is an $E$-convex and $\Phi$ is an $E$-convex of $u$ on $[0, \infty)$,
$$\Phi(E(t, 0)) = \lim_{u \to 0^+} \Phi(E(t,u)) = 0,$$
$$\lim_{u \to \infty} \Phi(E(t,u)) = \infty$$
and for each $u \in [0, \infty)$, $\Phi(E(t,u))$ is an $\mu$-measurable function of $t$ on $\Omega$.

**Remark 9.** Every Young function is an $E$-Young function if the map $E$ is taken as the identity map. But it is not every $E$-Young function a Young.

**Example 10.** Let $\Phi: \mathbb{R} \times [0, \infty) \to \mathbb{R}$ be defined by $\Phi(t,u) = e^{t+u} - 1$ and let $E: \mathbb{R} \times [0, \infty) \to \mathbb{R} \times [0, \infty)$ be defined by $E(t,u) = (u,u)$. Then, $\Phi$ is an $E$-Young function but it is not a Young function because for $\mu$-a.e. $t \in \mathbb{R}$, $\Phi(t, 0) = e^t - 1 \neq 0$.

**Example 11.** Let $\Phi: \mathbb{C} \times [0, \infty) \to \mathbb{R}$ be defined by
$$\Phi(t,u) = \begin{cases} t \ln(u), & u > 1 \\ 0, & 0 \leq u \leq 1 \end{cases}$$
and let $E: \mathbb{C} \times [0, \infty) \to \mathbb{C} \times [0, \infty)$ be defined by $E(t,u) = (-|t|, u)$. So, $\Phi$ is an $E$-Young function but it is not a Young function because, for $\mu$-a.e. $t \in \mathbb{C}$, $\Phi(t,u)$ is not convex because for $t \in (0, \infty)$, $\frac{\partial^2 \Phi}{\partial u^2} = -\frac{t}{u^2} < 0$.

**Definition 12.** Let $(\Omega, \Sigma, \mu)$ be a measure space. A function $\Phi: \Omega \times [0, \infty) \to \mathbb{R}$ is called an $E$-strong Young function if there exists a map $E: \Omega \times [0, \infty) \to \Omega \times [0, \infty)$ such that for $\mu$-a.e. $t \in \Omega$, $[0, \infty)$ is an $E$-convex and $\Phi$ is an $E$-convex $E$-continuous of $u$ on $[0, \infty)$, $\Phi(E(t, 0)) = 0 \Leftrightarrow u = 0$,
$$\lim_{u \to \infty} \Phi(E(t,u)) = \infty$$
and for each $u \in [0, \infty)$, $\Phi(E(t,u))$ is an $\mu$-measurable function of $t$ on $\Omega$.

**Remark 13.** Every strong Young function is an $E$-strong Young function if the map $E$ is taken as the identity map. But it is not every $E$-strong Young function a strong Young function.

**Example 14.** Let $\Phi: \mathbb{R} \times [0, \infty) \to \mathbb{R}$ be defined by $\Phi(t,u) = e^{u^t} - 1$ and let $E: \mathbb{R} \times [0, \infty) \to \mathbb{R} \times [0, \infty)$ be defined by $E(t,u) = (|t|, u)$. Then $\Phi$ is an $E$-strong Young function but it is not a strong Young function, where $\Phi(t,u) = e^{u^t} - 1$ is not convex because for $t \in (-\infty, 0)$ that $u^t$ is not convex.

**Example 15.** Let $\Phi: [0, \infty) \times [0, \infty) \to \mathbb{R}$ be defined by $\Phi(t,u) = \cosh(te^u) - 1$ and let $E: [0, \infty) \times [0, \infty) \to [0, \infty) \times [0, \infty)$ be defined by $E(t,u) = (u, 0)$. Then $\Phi$ is an $E$-strong Young function but it is not a strong Young function since for $\mu$-a.e. $t \in [0, \infty)$, $\Phi(t, 0) = \cosh(t) - 1 \neq 0$.

**Definition 16.** Let $(\Omega, \Sigma, \mu)$ be a measure space. A function $\Phi: \Omega \times [0, \infty) \to \mathbb{R}$ is called an $E$-Orlicz function if there exists a map $E: \Omega \times [0, \infty) \to \Omega \times [0, \infty)$ such that for $\mu$-a.e. $t \in \Omega$, $[0, \infty)$ is an $E$-convex and $\Phi$ is an $E$-convex of $u$ on $[0, \infty)$, $\Phi(E(t, 0)) = 0$, $\Phi(E(t,u)) > 0$ for any $u \in (0, \infty)$,
$$\lim_{u \to \infty} \Phi(E(t,u)) = \infty,$$
$\Phi$ is left $E$-continuous at
$$U_\Phi = \sup\{u > 0: \Phi(E(t,u)) < +\infty\}$$
and for each $u \in [0, \infty)$, $\Phi(E(t,u))$ is an $\mu$-measurable function of $t$ on $\Omega$.

**Remark 17.** Every Orlicz function is an $E$-Orlicz function if the map $E$ is taken as the identity map. But it is not every $E$-Orlicz function an Orlicz function.

**Example 18.** Let $\Phi: \mathbb{R} \times [0, \infty) \to \mathbb{R}$ be defined by $\Phi(t,u) = -t + u$ and let $E: \mathbb{R} \times [0, \infty) \to \mathbb{R} \times [0, \infty)$ be defined by $E(t,u) = (0, u^p), p \geq 1$. Then $\Phi$ is an $E$-Orlicz function but it is not an Orlicz function because for $\mu$-a.e. $t \in \mathbb{R}$, $\Phi(t, 0) = -t \neq 0$.

**Example 19.** Let $\Phi: \mathbb{R} \times [0, \infty) \to \mathbb{R}$ be defined by $\Phi(t,u) = t + u^{\frac{p}{(1-t)}}, p \geq 1$ and let $E: \mathbb{R} \times [0, \infty) \to \mathbb{R} \times [0, \infty)$ be defined by $E(t,u) = (0, u)$. Then $\Phi$ is an $E$-Orlicz function but it is not an Orlicz function because for $\mu$-a.e. $t \in \mathbb{R}$, $\Phi(t, 0) = t \neq 0$.

### III. ELEMENTARY PROPERTIES

A. Properties of $E$-$N$-Functions

**Theorem 20.** Let $\Phi_1, \Phi_2: \Omega \times [0, \infty) \to \mathbb{R}$ be $E$-$N$-functions with respect to $E: \Omega \times [0, \infty) \to \Omega \times [0, \infty)$. Then $\Phi_1 + \Phi_2$ and $c\Phi_1, c \geq 0$ are $E$-$N$-functions with respect to $E$.

**Theorem 21.** Let $\Phi\colon \Omega \times [0,\infty) \to \mathbb{R}$ be a linear $E$-$N$-function with respect to $E_1, E_2\colon \Omega \times [0,\infty) \to \Omega \times [0,\infty)$. Then $\Phi$ is an $E$-$N$-function with respect to $E_1 + E_2$ and $cE_1, c \geq 0$.

**Theorem 22.** Let $\Phi\colon \Omega \times [0,\infty) \to \mathbb{R}$ be a linear $E$-$N$-function with respect to $E_1, E_2\colon \Omega \times [0,\infty) \to \Omega \times [0,\infty)$. Then $\Phi$ is an $E$-$N$-function with respect to $E_1 \circ E_2$ and $E_2 \circ E_1$.

**Theorem 23.** Let $\Phi_i\colon \Omega \times [0,\infty) \to \mathbb{R}$ for $i = 1, \cdots, n$ be $E$-$N$-functions with respect to $E\colon \Omega \times [0,\infty) \to \Omega \times [0,\infty)$. Then $\Phi = \max_i \Phi_i$ is an $E$-$N$-function with respect to $E$.

**Theorem 24.** Let $\Phi\colon \Omega \times [0,\infty) \to \mathbb{R}$ be an $E$-$N$-function with respect to $E_i\colon \Omega \times [0,\infty) \to \Omega \times [0,\infty)$, $i = 1, \cdots, n$. Then $\Phi$ is an $E$-$N$-function with respect to $E_M = \max_i E_i$ and $E_m = \min_i E_i$.

**Theorem 25.** Let $(\Phi_n)_{n \in \mathbb{N}}$ be a sequence of continuous $E$-$N$-functions defined on a compact set $\Omega \times [0,\infty)$ with respect to $E\colon \Omega \times [0,\infty) \to \Omega \times [0,\infty)$ such that $(\Phi_n)_{n \in \mathbb{N}}$ converges uniformaly to a continuous function $\Phi\colon \Omega \times [0,\infty) \to \mathbb{R}$. Then $\Phi$ is an $E$-$N$-function with respect to $E$.

**Proof.** Assume that $(\Phi_n)_{n \in \mathbb{N}}$ is a sequence of continuous $E$-$N$-functions with respect to a map $E$ such that $\Phi_n \to \Phi$ uniformly on compact set $\Omega \times [0,\infty)$ and $\Phi$ is continuous on $\Omega \times [0,\infty)$. Then $\Phi_n(E) \to \Phi(E)$ uniformaly on $\Omega \times [0,\infty)$ and for $\mu$-a.e. $t \in \Omega$,
$$\Phi\big(E(t,u)\big) = \lim_{n \to \infty} \Phi_n\big(E(t,u)\big)$$
is even continuous convex of $u$ on $[0,\infty)$, $\Phi\big(E(t,u)\big) > 0$ for any $u \in (0,\infty)$,
$$\lim_{u \to 0} \frac{\Phi\big(E(t,u)\big)}{u} = \lim_{n \to \infty} \lim_{u \to 0} \frac{\Phi_n\big(E(t,u)\big)}{u} = 0,$$
$$\lim_{u \to \infty} \frac{\Phi\big(E(t,u)\big)}{u} = \lim_{n \to \infty} \lim_{u \to \infty} \frac{\Phi_n\big(E(t,u)\big)}{u} = \infty$$
and for each $u \in [0,\infty)$, $\Phi\big(E(t,u)\big)$ is an $\mu$-measurable function of $t$ on $\Omega$. ∎

**Theorem 26.** Let $\Phi$ be a continuous $E$-$N$-function defined on a compact set $\Omega \times [0,\infty)$ with respect to a sequence of maps $(E_n)_{n \in \mathbb{N}}$, $E_n\colon \Omega \times [0,\infty) \to \Omega \times [0,\infty)$, such that $(E_n)_{n \in \mathbb{N}}$ converges uniformaly to a map $E\colon \Omega \times [0,\infty) \to \Omega \times [0,\infty)$. Then $\Phi$ is an $E$-$N$-function with respect to $E$.

**Proof.** Suppose that $\Phi$ is a continuous $E$-$N$-function with respect to a sequence of maps $(E_n)_{n \in \mathbb{N}}$ such that $E_n \to E$ uniformaly on a compact set $\Omega \times [0,\infty)$. Then $\Phi(E_n) \to \Phi(E)$ uniformaly on $\Omega \times [0,\infty)$ and for $\mu$-a.e. $t \in \Omega$,
$$\Phi\big(E(t,u)\big) = \lim_{n \to \infty} \Phi\big(E_n(t,u)\big)$$
is even continuous convex of $u$ on $[0,\infty)$, $\Phi\big(E(t,u)\big) > 0$ for $u \in (0,\infty)$,
$$\lim_{u \to 0} \frac{\Phi\big(E(t,u)\big)}{u} = \lim_{n \to \infty} \lim_{u \to 0} \frac{\Phi\big(E_n(t,u)\big)}{u} = 0,$$
$$\lim_{u \to \infty} \frac{\Phi\big(E(t,u)\big)}{u} = \lim_{n \to \infty} \lim_{u \to \infty} \frac{\Phi\big(E_n(t,u)\big)}{u} = \infty$$
and for each $u \in [0,\infty)$, $\Phi\big(E(t,u)\big)$ is an $\mu$-measurable function of $t$ on $\Omega$. ∎

**Theorem 27.** Let $(\Phi_n)_{n \in \mathbb{N}}$ be a sequence of continuous $E$-$N$-functions defined on a compact set $\Omega \times [0,\infty)$ with respect to a sequence of continuous maps $(E_n)_{n \in \mathbb{N}}$, $E_n\colon \Omega \times [0,\infty) \to \Omega \times [0,\infty)$, such that $(\Phi_n)_{n \in \mathbb{N}}$ converges uniformaly to a continuous function $\Phi\colon \Omega \times [0,\infty) \to \mathbb{R}$ and $(E_n)_{n \in \mathbb{N}}$ converges uniformaly to a continuous map $E\colon \Omega \times [0,\infty) \to \Omega \times [0,\infty)$. Then $\Phi$ is an $E$-$N$-function with respect to $E$.

**Proof.** Assume that $(\Phi_n)_{n \in \mathbb{N}}$ is a sequence of continuous $E$-$N$-functions with respect to a sequence of continuous maps $(E_n)_{n \in \mathbb{N}}$ such that $\Phi_n \to \Phi$ uniformly and $E_n \to E$ uniformly on a compact set $\Omega \times [0,\infty)$ and $\Phi$ and $E$ are continuous on $\Omega \times [0,\infty)$. So $\Phi_n(E_n) \to \Phi(E)$ uniformaly on $\Omega \times [0,\infty)$ and for $\mu$-a.e. $t \in \Omega$, that
$$\Phi\big(E(t,u)\big) = \lim_{n \to \infty} \Phi_n\big(E_n(t,u)\big)$$
is even continuous convex of $u$ on $[0,\infty)$, $\Phi\big(E(t,u)\big) > 0, u \in (0,\infty)$,
$$\lim_{u \to 0} \frac{\Phi\big(E(t,u)\big)}{u} = \lim_{n \to \infty} \lim_{u \to 0} \frac{\Phi_n\big(E_n(t,u)\big)}{u} = 0,$$
$$\lim_{u \to \infty} \frac{\Phi\big(E(t,u)\big)}{u} = \lim_{n \to \infty} \lim_{u \to \infty} \frac{\Phi_n\big(E_n(t,u)\big)}{u} = \infty$$
and for each $u \in [0,\infty)$, $\Phi\big(E(t,u)\big)$ is an $\mu$-measurable function of $t$ on $\Omega$. ∎

B. Properties of $E$-Young Functions

**Theorem 28.** Let $\Phi_1, \Phi_2\colon \Omega \times [0,\infty) \to \mathbb{R}$ be $E$-Young functions with respect to $E\colon \Omega \times [0,\infty) \to \Omega \times [0,\infty)$. Then $\Phi_1 + \Phi_2$ and $c\Phi_1, c \geq 0$ are $E$-Young functions with respect to $E$.

**Theorem 29.** Let $\Phi\colon \Omega \times [0,\infty) \to \mathbb{R}$ be a linear $E$-Young function with respect to $E_1, E_2\colon \Omega \times [0,\infty) \to \Omega \times [0,\infty)$. Then $\Phi$ is an $E$-Young functions with respect to $E_1 + E_2$ and $cE_1, c \geq 0$.

**Theorem 30.** Let $\Phi\colon \Omega \times [0,\infty) \to \mathbb{R}$ be a linear $E$-Young function with respect to $E_1, E_2\colon \Omega \times [0,\infty) \to \Omega \times [0,\infty)$. Then $\Phi$ is an $E$-Young functions with respect to $E_1 \circ E_2$ and $E_2 \circ E_1$.

**Theorem 31.** Let $\Phi_i\colon \Omega \times [0,\infty) \to \mathbb{R}, i = 1, \ldots, n$ be $E$-Young functions with respect to $E\colon \Omega \times [0,\infty) \to \Omega \times [0,\infty)$. Then $\Phi = \max_i \Phi_i$ is an $E$-Young function with respect to $E$.

**Theorem 32.** Let $\Phi\colon \Omega \times [0,\infty) \to \mathbb{R}$ be an $E$-Young function with respect to $E_i\colon \Omega \times [0,\infty) \to \Omega \times [0,\infty)$, $i = 1, \cdots, n$. Then $\Phi$ is an $E$-Young function with respect to $E_M = \max_i E_i$ and $E_m = \min_i E_i$.

**Theorem 33.** Let $(\Phi_n)_{n \in \mathbb{N}}$ be a sequence of continuous $E$-Young functions defined on a compact set $\Omega \times [0,\infty)$ with respect to $E\colon \Omega \times [0,\infty) \to \Omega \times [0,\infty)$ such that $(\Phi_n)_{n \in \mathbb{N}}$ converges uniformaly to a continuous function $\Phi: \Omega \times [0, \infty) \to \mathbb{R}$. Then $\Phi$ is an $E$-Young function with respect to $E$.

**Proof.** Assume that $(\Phi_n)_{n \in \mathbb{N}}$ is a sequence of continuous $E$-Young functions with respect to a map $E$ such that $\Phi_n \to \Phi$ uniformly on a compact set $\Omega \times [0, \infty)$ and $\Phi$ is continuous on $\Omega \times [0, \infty)$. Then $\Phi_n(E) \to \Phi(E)$ uniformaly on $\Omega \times [0, \infty)$ and for $\mu$-a.e. $t \in \Omega$,
$$\Phi(E(t,u)) = \lim_{n \to \infty} \Phi_n(E(t,u))$$
is convex of $u$ on $[0, \infty)$,
$$\Phi(E(t,0)) = \lim_{u \to 0^+} \Phi(E(t,u))$$
$$= \lim_{n \to \infty} \lim_{u \to 0^+} \Phi_n(E(t,u)) = 0,$$
$$\lim_{u \to \infty} \Phi(E(t,u)) = \lim_{n \to \infty} \lim_{u \to \infty} \Phi_n(E(t,u)) = \infty$$
and for each $u \in [0, \infty)$, $\Phi(E(t,u))$ is an $\mu$-measurable function of $t$ on $\Omega$. ∎

**Theorem 34.** Let $\Phi: \Omega \times [0, \infty) \to \mathbb{R}$ be a continuous $E$-Young function defined on a compact set $\Omega \times [0, \infty)$ with respect to a sequence of maps $(E_n)_{n \in \mathbb{N}}$, $E_n: \Omega \times [0, \infty) \to \Omega \times [0, \infty)$, such that $(E_n)_{n \in \mathbb{N}}$ converges uniformaly to a map $E: \Omega \times [0, \infty) \to \Omega \times [0, \infty)$. Then $\Phi$ is an $E$-Young function with respect to $E$.

**Proof.** Suppose that $\Phi$ is a continuous $E$-Young function with respect to a sequence of maps $(E_n)_{n \in \mathbb{N}}$ such that $E_n \to E$ uniformly on a compact set $\Omega \times [0, \infty)$. Then $\Phi(E_n) \to \Phi(E)$ uniformly on $\Omega \times [0, \infty)$ and for $\mu$-a.e. $t \in \Omega$,
$$\Phi(E(t,u)) = \lim_{n \to \infty} \Phi(E_n(t,u))$$
is convex of $u$ on $[0, \infty)$,
$$\Phi(E(t,0)) = \lim_{u \to 0^+} \Phi(E(t,u))$$
$$= \lim_{n \to \infty} \lim_{u \to 0^+} \Phi(E_n(t,u)) = 0,$$
$$\lim_{u \to \infty} \Phi(E(t,u)) = \lim_{n \to \infty} \lim_{u \to \infty} \Phi(E_n(t,u)) = \infty$$
and for each $u \in [0, \infty)$, $\Phi(E(t,u))$ is an $\mu$-measurable function of $t$ on $\Omega$. ∎

**Theorem 35.** Let $(\Phi_n)_{n \in \mathbb{N}}$ be a sequence of continuous $E$-Young functions defined on a compact set $\Omega \times [0, \infty)$ with respect to a sequence of continuous maps $(E_n)_{n \in \mathbb{N}}$, $E_n: \Omega \times [0, \infty) \to \Omega \times [0, \infty)$, such that $(\Phi_n)_{n \in \mathbb{N}}$ converges uniformaly to a continuous function $\Phi: \Omega \times [0, \infty) \to \mathbb{R}$ and $(E_n)_{n \in \mathbb{N}}$ converges uniformaly to a continuous map $E: \Omega \times [0, \infty) \to \Omega \times [0, \infty)$. Then $\Phi$ is an $E$-Young function with respect to $E$.

**Proof.** Assume that $(\Phi_n)_{n \in \mathbb{N}}$ is a sequence of continuous $E$-Young functions with respect to a sequence of continuous maps $(E_n)_{n \in \mathbb{N}}$ such that $\Phi_n \to \Phi$ and $E_n \to E$ uniformly on a compact set $\Omega \times [0, \infty)$ and $\Phi$ and $E$ are continuous on $\Omega \times [0, \infty)$. Then $\Phi_n(E_n) \to \Phi(E)$ uniformaly on $\Omega \times [0, \infty)$ and for $\mu$-a.e. $t \in \Omega$,
$$\Phi(E(t,u)) = \lim_{n \to \infty} \Phi_n(E_n(t,u))$$
is convex of $u$ on $[0, \infty)$,
$$\Phi(E(t,0)) = \lim_{u \to 0^+} \Phi(E(t,u))$$
$$= \lim_{n \to \infty} \lim_{u \to 0^+} \Phi_n(E_n(t,u)) = 0,$$
$$\lim_{u \to \infty} \Phi(E(t,u)) = \lim_{n \to \infty} \lim_{u \to \infty} \Phi_n(E_n(t,u)) = \infty$$
and for each $u \in [0, \infty)$, $\Phi(E(t,u))$ is an $\mu$-measurable function of $t$ on $\Omega$. ∎

C. Properties of $E$-Strong Young Functions

**Theorem 36.** Let $\Phi_1, \Phi_2: \Omega \times [0, \infty) \to \mathbb{R}$ be $E$-strong Young functions with respect to $E: \Omega \times [0, \infty) \to \Omega \times [0, \infty)$. Then $\Phi_1 + \Phi_2$ and $c\Phi_1, c \geq 0$ are $E$-strong Young functions with respect to $E$.

**Theorem 37.** Let $\Phi: \Omega \times [0, \infty) \to \mathbb{R}$ be a linear $E$-strong Young function with respect to $E_1, E_2: \Omega \times [0, \infty) \to \Omega \times [0, \infty)$. Then $\Phi$ is an $E$-strong Young function with respect to $E_1 + E_2$ and $cE_1, c \geq 0$.

**Theorem 38.** Let $\Phi: \Omega \times [0, \infty) \to \mathbb{R}$ be a linear $E$-strong Young function with respect to $E_1, E_2: \Omega \times [0, \infty) \to \Omega \times [0, \infty)$. Then $\Phi$ is an $E$-strong Young function with respect to $E_1 \circ E_2$ and $E_2 \circ E_1$.

**Theorem 39.** Let $\Phi_i: \Omega \times [0, \infty) \to \mathbb{R}, i = 1, \cdots, n$ be $E$-strong Young functions with respect to $E: \Omega \times [0, \infty) \to \Omega \times [0, \infty)$. Then $\Phi = \max_i \Phi_i$ is an $E$-strong Young function with respect to $E$.

**Theorem 40.** Let $\Phi: \Omega \times [0, \infty) \to \mathbb{R}$ be an $E$-strong Young function with respect to $E_i: \Omega \times [0, \infty) \to \Omega \times [0, \infty), i = 1, \cdots, n$. Then $\Phi$ is an $E$-strong Young function with respect to $E_M = \max_i E_i$ and $E_m = \min_i E_i$.

**Theorem 41.** Let $(\Phi_n)_{n \in \mathbb{N}}$ be a sequence of continuous $E$-strong Young functions defined on a compact set $\Omega \times [0, \infty)$ with respect to $E: \Omega \times [0, \infty) \to \Omega \times [0, \infty)$ such that $(\Phi_n)_{n \in \mathbb{N}}$ converges uniformaly to a continuous function $\Phi: \Omega \times [0, \infty) \to \mathbb{R}$. Then $\Phi$ is an $E$-strong Young function with respect to $E$.

**Proof.** Assume that $(\Phi_n)_{n \in \mathbb{N}}$ is a sequence of continuous $E$-strong Young functions with respect to a map $E$ such that $\Phi_n \to \Phi$ uniformly on a compact set $\Omega \times [0, \infty)$ and $\Phi$ is continuous on $\Omega \times [0, \infty)$. Then $\Phi_n(E) \to \Phi(E)$ uniformaly on $\Omega \times [0, \infty)$ and for $\mu$-a.e. $t \in \Omega$,
$$\Phi(E(t,u)) = \lim_{n \to \infty} \Phi_n(E(t,u))$$
is convex, continuous of $u$ on $[0, \infty)$,
$$\Phi(E(t,0)) = \lim_{n \to \infty} \Phi_n(E(t,0)) = 0 \Leftrightarrow u = 0,$$
$$\lim_{u \to \infty} \Phi(E(t,u)) = \lim_{n \to \infty} \lim_{u \to \infty} \Phi_n(E(t,u)) = \infty$$
and for each $u \in [0, \infty)$, $\Phi(E(t,u))$ is an $\mu$-measurable function of $t$ on $\Omega$. ∎

**Theorem 42.** Let $\Phi: \Omega \times [0, \infty) \to \mathbb{R}$ be a continuous $E$-strong Young function defined on a compact set $\Omega \times [0, \infty)$ with respect to a sequence of maps $(E_n)_{n \in \mathbb{N}}$, $E_n: \Omega \times [0, \infty) \to \Omega \times [0, \infty)$ such that $(E_n)_{n \in \mathbb{N}}$ converges uniformaly to a map $E: \Omega \times [0, \infty) \to \Omega \times [0, \infty)$. Then $\Phi$ is an $E$-strong Young function with respect to $E$.

**Proof.** Suppose that $\Phi$ is a continuous $E$-strong Young function with respect to a sequence of maps $(E_n)_{n\in\mathbb{N}}$ such that $E_n \longrightarrow E$ uniformaly on a compact set $\Omega \times [0,\infty)$ and $E$ is continuous on $\Omega \times [0,\infty)$. Then $\Phi(E_n) \longrightarrow \Phi(E)$ uniformaly on $\Omega \times [0,\infty)$ and for $\mu$-a.e. $t \in \Omega$,
$$\Phi(E(t,u)) = \lim_{n\to\infty} \Phi(E_n(t,u))$$
is convex continuous of $u$ on $[0,\infty)$,
$$\Phi(E(t,0)) = \lim_{n\to\infty} \Phi(E_n(t,0)) = 0 \Leftrightarrow u = 0,$$
$$\lim_{u\to\infty} \Phi(E(t,u)) = \lim_{n\to\infty}\lim_{u\to\infty} \Phi(E_n(t,u)) = \infty$$
and for each $u \in [0,\infty), \Phi(E(t,u))$ is an $\mu$-measurable function of $t$ on $\Omega$. ∎

**Theorem 43.** Let $(\Phi_n)_{n\in\mathbb{N}}$ be a sequence of continuous $E$-strong Young functions defined on a compact set $\Omega \times [0,\infty)$ with respect to a sequence of continuous maps $(E_n)_{n\in\mathbb{N}}, E_n: \Omega \times [0,\infty) \longrightarrow \Omega \times [0,\infty)$ such that $(\Phi_n)_{n\in\mathbb{N}}$ converges uniformaly to a continuous function $\Phi: \Omega \times [0,\infty) \longrightarrow \mathbb{R}$ and $(E_n)_{n\in\mathbb{N}}$ converges uniformaly to a continuous map $E: \Omega \times [0,\infty) \longrightarrow \Omega \times [0,\infty)$. Then $\Phi$ is an $E$-strong Young function with respect to $E$.

**Proof.** Assume that $(\Phi_n)_{n\in\mathbb{N}}$ is a sequence of continuou $E$-strong Young functions with respect to a sequence of continuous maps $(E_n)_{n\in\mathbb{N}}$ such that $\Phi_n \longrightarrow \Phi$ and $E_n \longrightarrow E$ uniformly on a compact set $\Omega \times [0,\infty)$ and $\Phi$ and $E$ are continuous on $\Omega \times [0,\infty)$. So, $\Phi_n(E_n) \longrightarrow \Phi(E)$ uniformaly on $\Omega \times [0,\infty)$ and for $\mu$-a.e. $t \in \Omega$,
$$\Phi(E(t,u)) = \lim_{n\to\infty} \Phi_n(E_n(t,u))$$
is convex continuous of $u$ on $[0,\infty)$,
$$\Phi(E(t,0)) = \lim_{n\to\infty} \Phi_n(E_n(t,0)) = 0 \Leftrightarrow u = 0,$$
$$\lim_{u\to\infty} \Phi(E(t,u)) = \lim_{n\to\infty}\lim_{u\to\infty} \Phi_n(E_n(t,u)) = \infty$$
and for each $u \in [0,\infty), \Phi(E(t,u))$ is an $\mu$-measurable function of $t$ on $\Omega$. ∎

D. Properties of $E$-Orlicz Functions

**Theorem 44.** Let $\Phi_1, \Phi_2: \Omega \times [0,\infty) \longrightarrow \mathbb{R}$ be $E$-Orlicz functions with respect to $E: \Omega \times [0,\infty) \longrightarrow \Omega \times [0,\infty)$. Then $\Phi_1 + \Phi_2$ and $c\Phi_1, c \geq 0$ are $E$-Orlicz functions with respect to $E$.

**Theorem 45.** Let $\Phi: \Omega \times [0,\infty) \longrightarrow \mathbb{R}$ be a linear $E$-Orlicz function with respect to $E_1, E_2: \Omega \times [0,\infty) \longrightarrow \Omega \times [0,\infty)$. Then $\Phi$ is an $E$-Orlicz function with respect to $E_1 \circ E_2$ and $E_2 \circ E_1$.

**Theorem 46.** Let $\Phi: \Omega \times [0,\infty) \longrightarrow \mathbb{R}$ be a linear $E$-Orlicz function with respect to $E_1, E_2: \Omega \times [0,\infty) \longrightarrow \Omega \times [0,\infty)$. Then $\Phi$ is an $E$-Orlicz function with respect to $E_1 + E_2$ and $cE_1, c \geq 0$.

**Theorem 47.** Let $\Phi_i: \Omega \times [0,\infty) \longrightarrow \mathbb{R}, i = 1,\cdots,n$ be $E$-Orlicz functions with respect to $E: \Omega \times [0,\infty) \longrightarrow \Omega \times [0,\infty)$. Then $\Phi = \max_i \Phi_i$ is an $E$-Orlicz function with respect to $E$.

**Theorem 48.** Let $\Phi: \Omega \times [0,\infty) \longrightarrow \mathbb{R}$ be an $E$-Orlicz function with respect to $E_i: \Omega \times [0,\infty) \longrightarrow \Omega \times [0,\infty)$, $i = 1,\cdots,n$. Then $\Phi$ is an $E$-Orlicz function with respect to $E_M = \max_i E_i$ and $E_m = \min_i E_i$.

**Theorem 49.** Let $(\Phi_n)_{n\in\mathbb{N}}$ be a sequence of continuous $E$-Orlicz functions with respect to $E: \Omega \times [0,\infty) \longrightarrow \Omega \times [0,\infty)$ such that $(\Phi_n)_{n\in\mathbb{N}}$ converges uniformaly to a continuous function $\Phi: \Omega \times [0,\infty) \longrightarrow \mathbb{R}$. Then $\Phi$ is an $E$-Orlicz function withrespect to $E$.

**Proof.** Assume that $(\Phi_n)_{n\in\mathbb{N}}$ is a sequence of continuous $E$-Orlicz functions with respect to a map $E$ such that $\Phi_n \longrightarrow \Phi$ uniformly on a compact set $\Omega \times [0,\infty)$ and $\Phi$ is continuous on $\Omega \times [0,\infty)$. Then $\Phi_n(E) \longrightarrow \Phi(E)$ uniformaly on $\Omega \times [0,\infty)$ and for $\mu$-a.e. $t \in \Omega$,
$$\Phi(E(t,u)) = \lim_{n\to\infty} \Phi_n(E(t,u))$$
is convex of $u$ on $[0,\infty)$,
$$\Phi(E(t,0)) = \lim_{n\to\infty} \Phi_n(E(t,0)) = 0,$$
$$\lim_{u\to\infty} \Phi(E(t,u)) = \lim_{n\to\infty}\lim_{u\to\infty} \Phi_n(E(t,u)) = \infty,$$
$0 < \Phi(E(t,u)) < \infty$ for any $u \in (0,\infty)$, $\Phi(E(t,u))$ is left continuous at
$$U_\Phi = \sup\{u > 0: \Phi(E(t,u)) < +\infty\}.$$
and for each $u \in [0,\infty), \Phi(E(t,u))$ is an $\mu$-measurable function of $t$ on $\Omega$. ∎

**Theorem 50.** Let $\Phi: \Omega \times [0,\infty) \longrightarrow \mathbb{R}$ be a continuous $E$-Orlicz function defined on a compact set $\Omega \times [0,\infty)$ with respect to a sequence of maps $(E_n)_{n\in\mathbb{N}}, E_n: \Omega \times [0,\infty) \longrightarrow \Omega \times [0,\infty)$ such that $(E_n)_{n\in\mathbb{N}}$ converges uniformaly to a map $E: \Omega \times [0,\infty) \longrightarrow \Omega \times [0,\infty)$. Then $\Phi$ is an $E$-Orlicz function with respect to $E$.

**Proof.** Suppose that $\Phi$ is a continuous $E$-Orlicz function with respect to a sequence of continuous maps $(E_n)_{n\in\mathbb{N}}$ such that $E_n \longrightarrow E$ uniformaly on a compact set $\Omega \times [0,\infty)$ and $E$ is continuous on $\Omega \times [0,\infty)$. Then $\Phi(E_n) \longrightarrow \Phi(E)$ uniformaly on $\Omega \times [0,\infty)$ and for $\mu$-a.e. $t \in \Omega$,
$$\Phi(E(t,u)) = \lim_{n\to\infty} \Phi(E_n(t,u))$$
is convex of $u$ on $[0,\infty)$,
$$\Phi(E(t,0)) = \lim_{n\to\infty} \Phi(E_n(t,0)) = 0,$$
$$\lim_{u\to\infty} \Phi(E(t,u)) = \lim_{n\to\infty}\lim_{u\to\infty} \Phi(E_n(t,u)) = \infty,$$
$0 < \Phi(E(t,u)) < \infty$ for any $u \in (0,\infty)$ and $\Phi(E(t,u))$ is left continuous at
$$U_\Phi = \sup\{u > 0: \Phi(E(t,u)) < +\infty\}.$$
and for each $u \in [0,\infty), \Phi(E(t,u))$ is an $\mu$-measurable function of $t$ on $\Omega$. ∎

**Theorem 51.** Let $(\Phi_n)_{n\in\mathbb{N}}$ be a sequence of continuous $E$-Orlicz functions defined on a compact set $\Omega \times [0,\infty)$ with respect to a sequence of continuous maps $(E_n)_{n\in\mathbb{N}}, E_n: \Omega \times [0,\infty) \longrightarrow \Omega \times [0,\infty)$ such that $(\Phi_n)_{n\in\mathbb{N}}$ converges uniformaly to a continuous function $\Phi: \Omega \times [0,\infty) \longrightarrow \mathbb{R}$ and $(E_n)_{n\in\mathbb{N}}$ converges uniformaly to a continuous map $E: \Omega \times [0,\infty) \longrightarrow \Omega \times [0,\infty)$. Then $\Phi$ is an $E$-Orlicz function with respect to $E$.

**Proof.** Assume that $(\Phi_n)_{n\in\mathbb{N}}$ is a sequence of continuous $E$-Orlicz functions with respect to a sequence of continuous maps $(\Phi_n)_{n\in\mathbb{N}}$ such that $\Phi_n \to \Phi$ and $E_n \to E$ uniformly on a compact set $\Omega \times [0,\infty)$ and $\Phi$ and $E$ are continuous on $\Omega \times [0,\infty)$. Then $\Phi_n(E_n) \to \Phi(E)$ uniformaly on $\Omega \times [0,\infty)$ and for $\mu$-a.e. $t \in \Omega$,
$$\Phi(E(t,u)) = \lim_{n\to\infty} \Phi_n(E_n(t,u))$$
is convex of $u$ on $[0,\infty)$,
$$\Phi(E(t,0)) = \lim_{n\to\infty} \Phi_n(E_n(t,0)) = 0,$$
$$\lim_{u\to\infty} \Phi(E(t,u)) = \lim_{n\to\infty}\lim_{u\to\infty} \Phi_n(E_n(t,u)) = \infty,$$
$0 < \Phi(E(t,u)) < \infty$ for any $u \in (0,\infty)$, $\Phi(E)$ is left continuous at
$$U_\Phi = \sup\{u > 0: \Phi(E(t,u)) < +\infty\}.$$
and for each $u \in [0,\infty), \Phi(E(t,u))$ is an $\mu$-measurable function of $t$ on $\Omega$. ∎

## IV. RELATIONSHIPS BETWEEN $E$-CONVEX FUNCTIONS

In this section, we generalize the theorems in [9] to consider the relationships between $E$-$N$-functions, $E$-Young functions, $E$-strong Young functions and $E$-Orlicz functions.

**Theorem 52.** If $\Phi$ is an $E$-$N$-function, then $\Phi$ is an $E$-strong Young function.
**Proof.** Assume $\Phi: \Omega \times [0,\infty) \to \mathbb{R}$ is an $E$-$N$-function with a map $E: \Omega \times [0,\infty) \to \Omega \times [0,\infty)$. So, for $\mu$-a.e. $t \in \Omega, \Phi(E(t,u))$ is convex continuous of $u$ on $[0,\infty)$ satisfying
$$\forall \varepsilon > 0, \exists \delta > 0, 0 < u < \delta \Rightarrow \left|\frac{\Phi(E(t,u))}{u}\right| < \varepsilon$$
because
$$\lim_{u\to 0^+} \frac{\Phi(E(t,u))}{u} = 0.$$
Letting $\delta < 1$, we get
$$0 \leq |\Phi(E(t,u))| < \left|\frac{\Phi(E(t,u))}{\delta}\right| < \left|\frac{\Phi(E(t,u))}{u}\right| < \varepsilon.$$
By the squeeze theorem for functions, we get $\Phi(E(t,0)) = 0 \Leftrightarrow u = 0$ because $\Phi$ is continuous at $u = 0$ and $\Phi(E(t,u)) > 0$ for any $u \in (0,\infty)$. Moreover,
$$\forall M \in \mathbb{R}, \exists u_M > 0, u > u_M \Rightarrow \frac{\Phi(E(t,u))}{u} > M$$
because
$$\lim_{u\to\infty} \frac{\Phi(E(t,u))}{u} = \infty.$$
Taking $u_M > 1$, we have that
$$\Phi(E(t,u)) > Mu > Mu_M > M.$$
That is,
$$\lim_{u\to\infty} \Phi(E(t,u)) = \infty.$$
Furthermore, for each $u \in [0,\infty), \Phi(E(t,u))$ is an $\mu$-measurable function of $t$ on $\Omega$ which completes the proof. ∎

**Remark 53.** The converse of theorem 52 is not correct. That is, an $E$-strong Young function may not be an $E$-$N$-function.

**Example 54.** Let the function $\Phi: \mathbb{R} \times [0,\infty) \to \mathbb{R}$ be defined as $\Phi(t,u) = e^{u^t} - 1$ with the map $E: \mathbb{R} \times [0,\infty) \to \mathbb{R} \times [0,\infty)$ defined by $E(t,u) = (1,u)$. Then $\Phi$ is an $E$-strong Young function but it is not an $E$-$N$-function because for $\mu$-a.e. $t \in \mathbb{R}$,
$$\lim_{u\to 0} \frac{e^{u-1}}{u} = 1 \neq 0.$$

**Theorem 55.** If $\Phi$ is an $E$-strong Young function, then $\Phi$ is an $E$-Orlicz function.
**Proof.** Suppose that $\Phi: \Omega \times [0,\infty) \to \mathbb{R}$ is an $E$-strong Young function with a map $E: \Omega \times [0,\infty) \to \Omega \times [0,\infty)$. Then for $\mu$-a.e. $t \in \Omega, \Phi(E(t,u))$ is convex continuous of $u$ on $[0,\infty)$ satisfying $\Phi(E(t,0)) = 0$, $\Phi(E(t,u)) > 0$ for any $u \in (0,\infty)$ because $\Phi(E(t,0)) = 0 \Leftrightarrow u = 0$ and
$$\lim_{u\to\infty} \Phi(E(t,u)) = \infty$$
and $\Phi(E(t,u))$ is left continuous at $U_\Phi = +\infty$ because
$$\lim_{u\to\infty} \Phi(E(t,u)) = \infty.$$
Moreover, for each $u \in [0,\infty), \Phi(E(t,u))$ is an $\mu$-measurable function of $t$ on $\Omega$. Hence, $\Phi$ is an $E$-Orlicz function. ∎

**Remark 56.** The converse of theorem 55 is not correct. That is, an $E$-Orlicz function can not be an $E$-strong Young function.

**Example 57.** Let the function $\Phi: \mathbb{R} \times [0,\infty) \to \mathbb{R}$ be defined as
$$\Phi(t,u) = \begin{cases} u - |t|, 0 \leq u < 1 \\ u + |t| - 2, 1 \leq u \end{cases}$$
with the map $E: \mathbb{R} \times [0,\infty) \to \mathbb{R} \times [0,\infty)$ defined by $E(t,u) = (u,u)$. Then $\Phi$ is an $E$-Orlicz function but it not an $E$-strong Young function because, for $\mu$-a.e. $t \in \Omega, \Phi(E(t,1)) = 0$.

**Theorem 58.** If $\Phi$ is an $E$-Orlicz function, then $\Phi$ is an $E$-Young function.
**Proof.** Assume that $\Phi: \Omega \times [0,\infty) \to \mathbb{R}$ is an $E$-Orlicz function with a map $E: \Omega \times [0,\infty) \to \Omega \times [0,\infty)$. Then, for $\mu$-a.e. $t \in \Omega$, $\Phi(E(t,u))$ is convex of $u$ on $[0,\infty)$ satisfying $\Phi(E(t,0)) = 0$, $0 < \Phi(E(t,u))$, $u \in (0,\infty)$,
$$\lim_{u\to\infty} \Phi(E(t,u)) = \infty,$$
and $\Phi(E(t,u))$ is left continuous at $U_\Phi$. We only need to show that
$$\lim_{u\to 0^+} \Phi(E(t,u)) = 0.$$
In other words, we need to prove that
$$\forall \varepsilon > 0, \exists \delta_\varepsilon > 0, 0 < u < \delta_\varepsilon \Rightarrow 0 \leq \Phi(E(t,u)) < \varepsilon.$$
For arbitrary $\varepsilon > 0$, consider
$$a_\Phi = \inf\{u > 0: \Phi(E(t,u)) > 0\}.$$
If $a_\Phi > 0$, then $\Phi(E(t,u)) = 0$ for all $u \in (0, a_\Phi)$. Taking $\delta_\varepsilon = a_\Phi > 0$, then $\Phi(E(t,u)) = 0 < \varepsilon$ for all $0 < u < \delta_\varepsilon$. That is,

$$\lim_{u \to 0^+} \Phi(E(t,u)) = 0$$

If $a_\Phi = 0$, then $\Phi(E(t,u)) > 0$ for all $u > 0$ and there exists $0 < u_0 < \infty$ such that $0 < \Phi(E(t,u_0)) < \infty$. That is, for all $\varepsilon > 0, \exists u_\varepsilon \in (0,\infty)$ such that $0 < \Phi(E(t,u_\varepsilon)) < \infty$. If $\Phi(E(t,u_0)) < \varepsilon$, then $\Phi(E(t,u_\varepsilon)) < \infty$ for $u_\varepsilon = u_0$ and if $\Phi(E(t,u_0)) \geq \varepsilon$, then for $u_\varepsilon = \alpha u_0$ where $0 \leq \alpha = \frac{\varepsilon}{2\Phi(E(t,u_0))} < 1$ that
$$\Phi(E(t,u_\varepsilon)) = \Phi(E(t,\alpha u_0)) \leq \alpha \Phi(E(t,u_0))$$
$$\leq \frac{\varepsilon}{2} < \varepsilon$$

because $\Phi$ is $E$-convex of $u$ on $[0,\infty)$. Taking $\delta_\varepsilon = u_\varepsilon > 0$, we get, for $0 < u < \delta_\varepsilon$,
$$0 \leq \Phi(E(t,u)) \leq \Phi(E(t,\delta_\varepsilon)) = \Phi(E(t,u_\varepsilon)) < \varepsilon,$$
because $\Phi(E(t,u))$ is increasing of $u$ on $[0,\infty)$. Furthermore, for each $u \in [0,\infty), \Phi(E(t,u))$ is an $\mu$-measurable function of $t$ on $\Omega$. Hence, $\Phi$ is an $E$-Young function. ∎

**Remark 59.** The converse of theorem 58 is not correct. That is, it is not every $E$-Young function an $E$-Orlicz function.

**Example 60.** Let the function $\Phi: [0,\infty) \times [0,\infty) \to \mathbb{R}$ be defined as
$$\Phi(t,u) = \begin{cases} -\log(u + |t|^{1/p} + 1), & 0 \leq u < 1 \\ +\infty, & 1 \leq u \end{cases}$$
with the map $E: [0,\infty) \times [0,\infty) \to [0,\infty) \times [0,\infty)$ such that $E(t,u) = (u^p, u), p \geq 1$. Then $\Phi$ is an $E$-Young function but it is not an $E$-Orlicz function because $\Phi(E(t,u))$ is not left continuous at $U_\Phi = 1$ where
$$\lim_{u \to 1} \Phi(E(t,u)) = -\log(3) \neq +\infty$$
$$= \Phi(E(t,1)).$$

**Corollary 61.** $E$-$N$-function $\Rightarrow E$-strong Young function $\Rightarrow E$-Orlicz function $\Rightarrow E$-Young function.

**Corollary 62.** $E$-$N$-function $\nLeftarrow E$-strong Young function $\nLeftarrow E$-Orlicz function $\nLeftarrow E$-Young function.

## V. MAIN RESULTS

In this section, we are going to generate some of Orlicz spaces by $E$-Young functions and then establish their inclusion properties.

**Lemma 63.** Let $\Phi: \Omega \times [0,\infty) \to \mathbb{R}$ be an increasing $E$-Young function with respect to $E_1, E_2: \Omega \times [0,\infty) \to \Omega \times [0,\infty)$ such that, for $\mu$-a.e. $t \in \Omega$, $E_1(t,x) \leq E_2(t,x)$. Then, for $\mu$-a.e. $t \in \Omega$, $\Phi(E_1(t,x)) \leq \Phi(E_2(t,x))$.

**Lemma 64.** Let $\Phi_1, \Phi_2: \Omega \times [0,\infty) \to \mathbb{R}$ be $E$-Young functions with respect to $E: \Omega \times [0,\infty) \to \Omega \times [0,\infty)$ such that, for $\mu$-a.e. $t \in \Omega$, $\Phi_1(t,x) \leq \Phi_2(t,x)$. Then, for $\mu$-a.e. $t \in \Omega, \Phi_1(E(t,x)) \leq \Phi_2(E(t,x))$.

### A. E-Orlicz Spaces and Weak E-Orlicz Spaces

Let $\Phi: \Omega \times [0,\infty) \to \mathbb{R}$ be an $E$-Young function with respect to a map $E: \Omega \times [0,\infty) \to \Omega \times [0,\infty)$. The $E$-Orlicz space generated by $\Phi$ is defined by:
$$EL_{\Phi(E)}(\Omega, \Sigma, \mu) = \{f \in BS_\Omega: \|f\|_{\Phi(E)} < \infty\},$$
$$\|f\|_{\Phi(E)}$$
$$= \inf\left\{\lambda > 0: \int_\Omega \Phi\left(E\left(t, \frac{\|f(t)\|_{BS}}{\lambda}\right)\right) d\mu \leq 1\right\}$$
and the weak $E$-Orlicz space generated by $\Phi$ is
$$EL_{\Phi(E),weak}(\Omega, \Sigma, \mu) = \{f \in BS_\Omega: \|f\|_{\Phi(E),weak} < \infty\},$$
$$\|f\|_{\Phi(E),weak}$$
$$= \inf\left\{\lambda > 0: \sup_u \Phi(E(t,u)) \, m(\Omega, f/\lambda, u) \leq 1\right\}$$
where $BS_\Omega$ is the set of all $\mu$-measurable functions $f$ from $\Omega$ to $BS$ such that $(BS, \|\cdot\|_{BS})$ is a Banach space and
$$m(\Omega, f, u) = \mu\{t \in \Omega: \|f(t)\|_{BS} > u\}.$$

**Example 65.** We have seen from example 10 that $\Phi(t,u) = e^{t+u} - 1$ is an $E$-Young function with respect to the map $E(t,u) = (u,u)$. Then the $E$-Orlicz space and the weak $E$-Orlicz space generated by $\Phi(E(t,u)) = e^{2u} - 1$ are equipped with the norm
$$\|f\|_{\Phi(E)} = \inf\left\{\lambda > 0: \int_\Omega \left(\exp\left(\frac{2\|f(t)\|_{BS}}{\lambda}\right) - 1\right) d\mu \leq 1\right\},$$
for all $f \in EL_{\Phi(E)}(\Omega, \Sigma, \mu)$ and
$$\|f\|_{\Phi(E),weak}$$
$$= \inf\left\{\lambda > 0: \sup_u (e^{2u} - 1) \, m(\Omega, f/\lambda, u) \leq 1\right\}$$
for all $f \in EL_{\Phi(E),weak}(\Omega, \Sigma, \mu)$.
If $\Phi_p(E(t,u)) = u^p, p \geq 1$, we get
$$EL_p(\Omega, \Sigma, \mu) = EL_{\Phi_p(E)}(\Omega, \Sigma, \mu)$$
$$= \{f \in X_\Omega: \|f\|_p < \infty\},$$
$$\|f\|_p = \inf\left\{\lambda > 0: \int_\Omega |f/\lambda|^p d\mu \leq 1\right\}$$
$$= \left(\int_\Omega |f|^p d\mu\right)^{1/p}$$
for all $f \in EL_p(\Omega, \Sigma, \mu)$ and
$$EL_{p,weak}(\Omega, \Sigma, \mu) = EL_{\Phi_p,weak}(\Omega, \Sigma, \mu)$$
$$= \{f \in X_\Omega: \|f\|_{p,weak} < \infty\},$$
$$\|f\|_{p,weak} = \inf\left\{\lambda > 0: \sup_u u^p \, m(\Omega, f/\lambda, u) \leq 1\right\}.$$

**Example 66.** Let $\Phi: \mathbb{C} \times [0,\infty) \to \mathbb{R}$ be defined as
$$\Phi(t,u) = \begin{cases} t\ln(u), & u > 1 \\ 0, & 0 \leq u \leq 1 \end{cases}$$
with respect to $E: \mathbb{C} \times [0,\infty) \to \mathbb{C} \times [0,\infty)$ such that
$$E(t,u) = \begin{cases} (1, e^{u^p}), & 1 \leq p, \\ (1,0), & 1 < u, p = +\infty, \\ (0,0), & 0 \leq u \leq 1, p = +\infty. \end{cases}$$
Then, for $\mu$-a.e. $t \in \mathbb{C}$, that

$$\Phi(E(t,u)) = \begin{cases} u^p, & 1 \le p \\ +\infty, 1 < u, p = +\infty \\ 0, 0 \le u \le 1, p = +\infty \end{cases}$$

is an $E$-Young function and the obtained spaces are $EL_p(\Omega,\Sigma,\mu)$ and $EL_{p,weak}(\Omega,\Sigma,\mu)$ for $1 \le p \le \infty$.

**Theorem 67.** If $\Phi: \Omega \times [0,\infty) \to \mathbb{R}$ is an increasing $E$-Young function with respect to $E_1, E_2: \Omega \times [0,\infty) \to \Omega \times [0,\infty)$ such that, for $\mu$-a.e. $t \in \Omega$, $E_1(t,x) \le E_2(t,x)$. Then
$$EL_{\Phi(E_2)}(\Omega,\Sigma,\mu) \subseteq EL_{\Phi(E_1)}(\Omega,\Sigma,\mu)$$
and
$$EL_{\Phi(E_2),weak}(\Omega,\Sigma,\mu) \subseteq EL_{\Phi(E_1),weak}(\Omega,\Sigma,\mu).$$

**Theorem 68.** If $\Phi_1, \Phi_2: \Omega \times [0,\infty) \to \mathbb{R}$ are $E$-Young functions with respect to $E: \Omega \times [0,\infty) \to \Omega \times [0,\infty)$ such that, for $\mu$-a.e. $t \in \Omega$, $\Phi_1(E(t,x)) \le \Phi_2(E(t,x))$. Then
$$EL_{\Phi_2(E)}(\Omega,\Sigma,\mu) \subseteq EL_{\Phi_1(E)}(\Omega,\Sigma,\mu)$$
and
$$EL_{\Phi_2(E),weak}(\Omega,\Sigma,\mu) \subseteq EL_{\Phi_1(E),weak}(\Omega,\Sigma,\mu).$$

**Theorem 69.** If $\Phi: \Omega \times [0,\infty) \to \mathbb{R}$ is an increasing $E$-Young function with respect to $E_1, E_2: \Omega \times [0,\infty) \to \Omega \times [0,\infty)$ such that, for $\mu$-a.e. $t \in \Omega$, $E_1(t,x) \le E_2(t,x)$. Then
$$EL_{\Phi(E_2)}(\Omega,\Sigma,\mu) \subseteq EL_{\Phi(E_1),weak}(\Omega,\Sigma,\mu)$$
and if $\Omega \times [0,\infty)$ is compact, then
$$EL_{\Phi(E_2),weak}(\Omega,\Sigma,\mu) \subseteq EL_{\Phi(E_1)}(\Omega,\Sigma,\mu).$$

**Proof.** Let $f \in EL_{\Phi(E_2)}(\Omega,\Sigma,\mu)$ and let $\Phi$ be an increasing $E$-Young function. Then, by Lemma 59, we have
$$\Phi(E_1(t,u))m(\Omega,f/\lambda,u) \le \Phi(E_2(t,u))m(\Omega,f/\lambda,u)$$
$$\le \int_{\{t\in\Omega: \frac{\|f(t)\|_{BS}}{\lambda} > u\}} \Phi\left(E_2\left(t, \frac{\|f(t)\|_{BS}}{\lambda}\right)\right) d\mu$$
$$\le \int_\Omega \Phi\left(E_2\left(t, \frac{\|f(t)\|_{BS}}{\lambda}\right)\right) d\mu \le 1.$$

Since $u$ is arbitrary, we have
$$\sup_u \Phi(E_1(t,u)) m(\Omega,f/\lambda,u) \le 1$$
and $f \in EL_{\Phi(E_1),weak}(\Omega,\Sigma,\mu)$ with
$$\|f\|_{\Phi(E_1),weak} \le \|f\|_{\Phi(E_2)}.$$

Let $f \in EL_{\Phi(E_2),weak}(\Omega,\Sigma,\mu)$ and assume that $\Omega \times [0,\infty)$ is compact. Then
$$\int_\Omega \Phi\left(E_1\left(t, \frac{\|f(t)\|_{BS}}{\lambda}\right)\right) d\mu$$
$$= \sup_u \Phi(E_1(t,u)) m(\Omega,f/\lambda,u)$$
$$\le \sup_u \Phi(E_2(t,u)) m(\Omega,f/\lambda,u) \le 1.$$

That is, $f \in EL_{\Phi(E_1)}(\Omega,\Sigma,\mu)$ with
$$\|f\|_{\Phi(E_1)} \le \|f\|_{\Phi(E_2),weak}. \blacksquare$$

**Theorem 70.** If $\Phi_1, \Phi_2: \Omega \times [0,\infty) \to \mathbb{R}$ are $E$-Young functions with respect to $E: \Omega \times [0,\infty) \to \Omega \times [0,\infty)$ such that, for $\mu$-a.e. $t \in \Omega$, $\Phi_1(E(t,x)) \le \Phi_2(E(t,x))$. Then
$$EL_{\Phi_2(E)}(\Omega,\Sigma,\mu) \subseteq EL_{\Phi_1(E),weak}(\Omega,\Sigma,\mu)$$
and if $\Omega \times [0,\infty)$ is compact, then
$$EL_{\Phi_2(E),weak}(\Omega,\Sigma,\mu) \subseteq EL_{\Phi_1(E)}(\Omega,\Sigma,\mu).$$

**Proof.** Let $f \in EL_{\Phi_2(E)}(\Omega,\Sigma,\mu)$. Then
$$\Phi_1(E(t,u))m(\Omega,f/\lambda,u)$$
$$\le \Phi_2(E(t,u))m(\Omega,f/\lambda,u)$$
$$\le \int_{\{t\in\Omega: \frac{\|f(t)\|_{BS}}{\lambda} > u\}} \Phi_2\left(E\left(t, \frac{\|f(t)\|_{BS}}{\lambda}\right)\right) d\mu$$
$$\le \int_\Omega \Phi_2\left(E\left(t, \frac{\|f(t)\|_{BS}}{\lambda}\right)\right) d\mu \le 1.$$

Since $u$ is arbitrary, we have
$$\sup_u \Phi_1(E(t,u)) m(\Omega,f/\lambda,u) \le 1$$
and $f \in EL_{\Phi_1(E),weak}(\Omega,\Sigma,\mu)$ with
$$\|f\|_{\Phi_1(E),weak} \le \|f\|_{\Phi_2(E)}.$$

Let $f \in EL_{\Phi_2(E),weak}(\Omega,\Sigma,\mu)$ and assume that $\Omega \times [0,\infty)$ is compact. Then
$$\int_\Omega \Phi_1\left(E\left(t, \frac{\|f(t)\|_{BS}}{\lambda}\right)\right) d\mu$$
$$= \sup_u \Phi_1(E(t,u)) m(\Omega,f/\lambda,u)$$
$$\le \sup_u \Phi_2(E(t,u)) m(\Omega,f/\lambda,u) \le 1.$$

That is, $f \in EL_{\Phi_1(E)}(\Omega,\Sigma,\mu)$ with
$$\|f\|_{\Phi_1(E)} \le \|f\|_{\Phi_2(E),weak}. \blacksquare$$

**B. $E$-Orlicz-Sobolev Space and Weak $E$-Orlicz-Sobolev Space**

Let $\Phi: \Omega \times [0,\infty) \to \mathbb{R}$ be an $E$-Young function with respect to $E: \Omega \times [0,\infty) \to \Omega \times [0,\infty)$. The $E$-Orlicz-Sobolev space $EW^k L_{\Phi(E)}(\Omega,\Sigma,\mu)$ generated by $\Phi(E)$ is
$$EW^k L_{\Phi(E)}(\Omega,\Sigma,\mu) = \{f \in EL_{\Phi(E)}(\Omega,\Sigma,\mu):$$
$$D^\alpha f \in EL_{\Phi(E)}(\Omega,\Sigma,\mu),$$
$$\forall |\alpha| \le k\},$$
$$\|f\|_{k,\Phi(E)} = \sum_{|\alpha| \le k} \|D^\alpha f\|_{\Phi(E)}$$

for all $f \in EW^k L_{\Phi(E)}(\Omega,\Sigma,\mu)$ and the weak $E$-Orlicz-Sobolev space is
$$EW^k L_{\Phi(E),weak}(\Omega,\Sigma,\mu) = \{f \in EL_{\Phi(E),weak}(\Omega,\Sigma,\mu):$$
$$D^\alpha f \in EL_{\Phi(E),weak}(\Omega,\Sigma,\mu), \forall |\alpha| \le k\},$$
$$\|f\|_{k,\Phi(E),weak} = \sum_{|\alpha| \le k} \|D^\alpha f\|_{\Phi(E),weak}$$

for all $f \in EW^k L_{\Phi(E),weak}(\Omega,\Sigma,\mu)$.

If $\Phi_p(E(t,u)) = u^p, p \ge 1$, we get the $E$-Sobolev space
$$EW^k L_{\Phi_p(E)}(\Omega,\Sigma,\mu) = EW^{k,p}(\Omega,\Sigma,\mu)$$
$$= \{f \in EL_p(\Omega,\Sigma,\mu): D^\alpha f \in$$
$$EL_p(\Omega,\Sigma,\mu), \forall |\alpha| \le k\}$$

equipped with the norm
$$\|f\|_{k,p} = \left(\sum_{|\alpha| \le k} \|D^\alpha f\|_p\right)^{1/p}$$

for all $f \in EL_p(\Omega,\Sigma,\mu)$ and the weak $E$-Sobolev space

$$EW^kL_{\Phi_p(E),weak}(\Omega,\Sigma,\mu) = EW^{k,p,weak}(\Omega,\Sigma,\mu)$$
$$= \{f \in EL_{p,weak}(\Omega,\Sigma,\mu): D^\alpha f$$
$$\in EL_{p,weak}(\Omega,\Sigma,\mu), \forall\, |\alpha| \leq k\},$$
$$\|f\|_{k,p,weak} = \sum_{|\alpha|\leq k} \|D^\alpha f\|_{p,weak}$$
for all $f \in EW^k L_{p,weak}(\Omega,\Sigma,\mu)$.

**Theorem 71.** If $\Phi:\Omega \times [0,\infty) \to \mathbb{R}$ is an increasing $E$-Young function with respect to $E_1, E_2: \Omega \times [0,\infty) \to \Omega \times [0,\infty)$ such that, for $\mu$-a.e. $t \in \Omega$, $E_1(t,x) \leq E_2(t,x)$. Then
$$EW^k L_{\Phi(E_2)}(\Omega,\Sigma,\mu) \subseteq EW^k L_{\Phi(E_1)}(\Omega,\Sigma,\mu)$$
and
$$EW^k L_{\Phi(E_2),weak}(\Omega,\Sigma,\mu) \subseteq EW^k L_{\Phi(E_1),weak}(\Omega,\Sigma,\mu).$$

**Theorem 72.** If $\Phi_1,\Phi_2:\Omega \times [0,\infty) \to \mathbb{R}$ are $E$-Young functions with respect to $E:\Omega \times [0,\infty) \to \Omega \times [0,\infty)$ such that, for $\mu$-a.e. $t \in \Omega$, $\Phi_1(E(t,x)) \leq \Phi_2(E(t,x))$. Then
$$EW^k L_{\Phi_2(E)}(\Omega,\Sigma,\mu) \subseteq EW^k L_{\Phi_1(E)}(\Omega,\Sigma,\mu)$$
and
$$EW^k L_{\Phi_2(E),weak}(\Omega,\Sigma,\mu) \subseteq EW^k L_{\Phi_1(E),weak}(\Omega,\Sigma,\mu).$$

**Theorem 73.** If $\Phi:\Omega \times [0,\infty) \to \mathbb{R}$ is an increasing $E$-Young function with respect to $E_1, E_2:\Omega \times [0,\infty) \to \Omega \times [0,\infty)$ such that, for $\mu$-a.e. $t \in \Omega$, $E_1(t,x) \leq E_2(t,x)$. Then
$$EW^k L_{\Phi(E_2)}(\Omega,\Sigma,\mu) \subseteq EW^k L_{\Phi(E_1),weak}(\Omega,\Sigma,\mu)$$
and if $\Omega \times [0,\infty)$ is a compact set, then
$$EW^k L_{\Phi(E_2),weak}(\Omega,\Sigma,\mu)$$
$$\subseteq EW^k L_{\Phi(E_1)}(\Omega,\Sigma,\mu)(\Omega,\Sigma,\mu).$$

**Theorem 74.** If $\Phi_1,\Phi_2:\Omega \times [0,\infty) \to \mathbb{R}$ are $E$-Young functions with respect to $E:\Omega \times [0,\infty) \to \Omega \times [0,\infty)$ such that, for $\mu$-a.e. $t \in \Omega$, $\Phi_1(E(t,x)) \leq \Phi_2(E(t,x))$. Then
$$EW^k L_{\Phi_2(E)}(\Omega,\Sigma,\mu) \subseteq EW^k L_{\Phi_1(E),weak}(\Omega,\Sigma,\mu)$$
and if $\Omega \times [0,\infty)$ is a compact set, then
$$EW^k L_{\Phi_2(E),weak}(\Omega,\Sigma,\mu) \subseteq EW^k L_{\Phi_1(E)}(\Omega,\Sigma,\mu).$$

### C. E-Orlicz-Morrey Space and Weak E-Orlicz-Morrey Space

Let $\Phi:\Omega \times [0,\infty) \to \mathbb{R}$ be an $E$-convex function with respect to $E:\Omega \times [0,\infty) \to \Omega \times [0,\infty)$ and let $\phi:(0,\infty) \to (0,\infty)$ be a function such that $\phi(r)$ is almost decreasing and $\phi(r)r^n$ is almost increasing and let $B$ denote the ball $B(a,r) = \{t \in \Omega: |t-a| < r\}$. The $E$-Orlicz-Morrey space is
$$EL_{\Phi(E),\phi}(\Omega,\Sigma,\mu) = \{f \in X_\Omega: \|f\|_{\Phi(E),\phi} < \infty\},$$
$$\|f\|_{\Phi(E),\phi}$$
$$= \sup_B \inf\left\{\lambda > 0: \frac{1}{|B|\phi(r)}\int_B \Phi\left(E\left(t,\frac{\|f(t)\|_{BS}}{\lambda}\right)\right)d\mu \leq 1\right\},$$
and the weak $E$-Orlicz-Morrey space is
$$EL_{\Phi(E),\phi,weak}(\Omega,\Sigma,\mu)$$
$$= \{f \in X_\Omega: \|f\|_{\Phi(E),\phi,weak} < \infty\},$$
$$\|f\|_{\Phi(E),\phi,weak}$$
$$= \sup_B \inf\left\{\lambda > 0: \sup_u \frac{\Phi(E(t,u))m(B,f/\lambda,u)}{|B|\phi(r)} \leq 1\right\}.$$
If $\Phi_p(E(t,u)) = u^p, p \geq 1$, then
$$EL_{\Phi_p(E),\phi}(\Omega,\Sigma,\mu) = EL_{p,\phi}(\Omega,\Sigma,\mu)$$
$$= \{f \in X_\Omega: \|f\|_{p,\phi} < \infty\},$$
$$\|f\|_{p,\phi} = \sup_B \left(\frac{1}{|B|\phi(r)}\int_B \|f(t)\|_{BS}^p d\mu\right)^{1/p},$$
$$EL_{\Phi_p(E),\phi,weak}(\Omega,\Sigma,\mu) = EL_{p,\phi,weak}(\Omega,\Sigma,\mu)$$
$$= \{f \in X_\Omega: \|f\|_{p,\phi,weak} < \infty\},$$
$$\|f\|_{p,\phi,weak} = \sup_B \sup_u \frac{u^p m(B,f,u)}{|B|\phi(r)}.$$
If $\phi(r) = r^{-n}$, we get
$$EL_{\Phi(E),\phi}(\Omega,\Sigma,\mu) = EL_{\Phi(E)}(\Omega,\Sigma,\mu),$$
$$EL_{\Phi(E),\phi,weak}(\Omega,\Sigma,\mu) = EL_{\Phi(E),weak}(\Omega,\Sigma,\mu).$$
If $\Phi_p(E(t,u)) = u^p, p \geq 1$ and $\phi(r) = r^{\lambda-n}$, we get the Morrey space
$$EL_{\Phi_p(E),\phi}(\Omega,\Sigma,\mu) = EL_{p,\lambda}(\Omega,\Sigma,\mu)$$
$$= \{f \in X_\Omega: \|f\|_{p,\lambda} < \infty\},$$
$$\|f\|_{p,\lambda} = \sup_B \left(\frac{1}{r^\lambda}\int_B \|f(t)\|_{BS}^p d\mu\right)^{1/p}$$
and the weak Morrey space is
$$EL_{\Phi_p(E),\phi,weak}(\Omega,\Sigma,\mu) = EL_{p,\lambda,weak}(\Omega,\Sigma,\mu)$$
$$= \{f \in X_\Omega: \|f\|_{p,\lambda,weak} < \infty\},$$
$$\|f\|_{p,\lambda,weak} = \sup_B \sup_u \frac{u^p m(B,f,u)}{r^\lambda}.$$

**Theorem 75.** If $\Phi:\Omega \times [0,\infty) \to \mathbb{R}$ is an increasing $E$-Young function with respect to $E_1, E_2:\Omega \times [0,\infty) \to \Omega \times [0,\infty)$ such that, for $\mu$-a.e. $t \in \Omega$, $E_1(t,x) \leq E_2(t,x)$. Then
$$EL_{\Phi(E_2),\phi}(\Omega,\Sigma,\mu) \subseteq EL_{\Phi(E_1),\phi}(\Omega,\Sigma,\mu)$$
and
$$EL_{\Phi(E_2),\phi,weak}(\Omega,\Sigma,\mu) \subseteq EL_{\Phi(E_1),\phi,weak}(\Omega,\Sigma,\mu).$$

**Theorem 76.** If $\Phi_1,\Phi_2:\Omega \times [0,\infty) \to \mathbb{R}$ are $E$-Young function with respect to $E:\Omega \times [0,\infty) \to \Omega \times [0,\infty)$ such that, for $\mu$-a.e. $t \in \Omega$, $\Phi_1(E(t,x)) \leq \Phi_2(E(t,x))$. Then
$$EL_{\Phi_2(E),\phi}(\Omega,\Sigma,\mu) \subseteq EL_{\Phi_1(E),\phi}(\Omega,\Sigma,\mu)$$
and
$$EL_{\Phi_2(E),\phi,weak}(\Omega,\Sigma,\mu) \subseteq EL_{\Phi_1(E),\phi,weak}(\Omega,\Sigma,\mu).$$

**Theorem 77.** If $\Phi:\Omega \times [0,\infty) \to \mathbb{R}$ is an increasing $E$-Young function with respect to $E_1, E_2:\Omega \times [0,\infty) \to \Omega \times [0,\infty)$ such that, for $\mu$-a.e. $t \in \Omega$, $E_1(t,x) \leq E_2(t,x)$. Then
$$EL_{\Phi(E_2),\phi}(\Omega,\Sigma,\mu) \subseteq EL_{\Phi(E_1),\phi,weak}(\Omega,\Sigma,\mu)$$
and if $\Omega \times [0,\infty)$ is a compact set, then
$$EL_{\Phi(E_2),\phi,weak}(\Omega,\Sigma,\mu) \subseteq EL_{\Phi(E_1),\phi}(\Omega,\Sigma,\mu).$$

**Proof.** Let $f \in EL_{\Phi(E_2),\phi}(\Omega,\Sigma,\mu)$ and let $\Phi$ be an increasing $E$-Young function. By Lemma 59, we have
$$\frac{\Phi(E_1(t,u))m(B,f/\lambda,u)}{|B|\phi(r)} \leq \frac{\Phi(E_2(t,u))m(B,f/\lambda,u)}{|B|\phi(r)}$$

$$\leq \frac{1}{|B|\phi(r)}\int_B \Phi\left(E_2\left(t,\frac{\|f(t)\|_{BS}}{\lambda}\right)\right)d\mu \leq 1.$$

Since $u$ is arbitrary, then
$$\sup_{u>0}\frac{\Phi(E_1(t,u))m(B,f/\lambda,u)}{|B|\phi(r)} \leq 1$$
and $f \in EL_{\Phi(E_1),\phi,weak}(\Omega,\Sigma,\mu)$ with
$$\|f\|_{\Phi(E_1),\phi,weak} \leq \|f\|_{\Phi(E_2),\phi}.$$
Let $f \in EL_{\Phi(E_2),\phi,weak}(\Omega,\Sigma,\mu)$ and $\Omega \times [0,\infty)$ be a compact set. Then
$$\frac{1}{|B|\phi(r)}\int_B \Phi\left(E_1\left(t,\frac{\|f(t)\|_{BS}}{\lambda}\right)\right)d\mu$$
$$= \sup_{u>0}\frac{\Phi(E_1(t,u))m(B,f/\lambda,u)}{|B|\phi(r)}$$
$$\leq \sup_{u>0}\frac{\Phi(E_2(t,u))m(B,f/\lambda,u)}{|B|\phi(r)} \leq 1.$$
So, $f \in EL_{\Phi(E_1),\phi}(\Omega,\Sigma,\mu)$ with
$$\|f\|_{\Phi(E_1)} \leq \|f\|_{\Phi(E_2),weak}.\blacksquare$$

**Theorem 78.** If $\Phi_1,\Phi_2:\Omega\times[0,\infty)\to\mathbb{R}$ are $E$-Young functions with respect to $E:\Omega\times[0,\infty)\to\Omega\times[0,\infty)$ such that, for $\mu$-a.e. $t\in\Omega$, $\Phi_1(E(t,x))\leq \Phi_2(E(t,x))$. Then
$$EL_{\Phi_2(E),\phi}(\Omega,\Sigma,\mu) \subseteq EL_{\Phi_1(E),\phi,weak}(\Omega,\Sigma,\mu)$$
and if $\Omega\times[0,\infty)$ is compact set, then
$$EL_{\Phi_2(E),\phi,weak}(\Omega,\Sigma,\mu) \subseteq EL_{\Phi_1(E),\phi}(\Omega,\Sigma,\mu).$$
**Proof.** Let $f \in EL_{\Phi_2(E),\phi}(\Omega,\Sigma,\mu)$. Then
$$\frac{\Phi_1(E(t,u))m(B,f/\lambda,u)}{|B|\phi(r)} \leq \frac{\Phi_2(E(t,u))m(B,f/\lambda,u)}{|B|\phi(r)}$$
$$\leq \frac{1}{|B|\phi(r)}\int_B \Phi_2\left(E\left(t,\frac{\|f(t)\|_{BS}}{\lambda}\right)\right)d\mu \leq 1.$$
Since $u$ is arbitrary, we have
$$\sup_{u>0}\frac{\Phi_1(E(t,u))m(B,f/\lambda,u)}{|B|\phi(r)} \leq 1$$
and $f \in EL_{\Phi_1(E),\phi,weak}(\Omega,\Sigma,\mu)$ with
$$\|f\|_{\Phi_1(E),\phi,weak} \leq \|f\|_{\Phi_2(E),\phi}.$$
Let $f \in EL_{\Phi_2(E),\phi,weak}(\Omega,\Sigma,\mu)$ and $\Omega\times[0,\infty)$ be a compact set. Then
$$\frac{1}{|B|\phi(r)}\int_B \Phi_1\left(E\left(t,\frac{\|f(t)\|_{BS}}{\lambda}\right)\right)d\mu$$
$$= \sup_{u>0}\frac{\Phi_1(E(t,u))m(B,f/\lambda,u)}{|B|\phi(r)}$$
$$\leq \sup_{u>0}\frac{\Phi_2(E(t,u))m(B,f/\lambda,u)}{|B|\phi(r)} \leq 1.$$
So, $f \in EL_{\Phi_1(E),\phi}(\Omega,\Sigma,\mu)$ with
$$\|f\|_{\Phi_1(E)} \leq \|f\|_{\Phi_2(E),weak}.\blacksquare$$

*D.  E-Orlicz-Lorentz Spaces*

Let $\Phi:[0,\infty)\times[0,\infty)\to\mathbb{R}$ be an $E$-convex function with respect to $E:(0,\infty)\times[0,\infty)\to(0,\infty)\times[0,\infty)$ and let $\omega:[0,\infty)\to[0,\infty)$ be a weight function and $W(t)=\int_0^t \omega(s)ds$. The $E$-Orlicz-Lorentz space is:
$$\Lambda_{\omega,\Phi(E)} = \{f\in X_{(0,\infty)}: \|f\|_{\omega,\Phi(E)} < \infty\},$$
$$\|f\|_{\omega,\Phi(E)}$$
$$= \inf\left\{\lambda>0: \int_0^\infty \Phi(E(t,f^*(t)/\lambda))W(t)d\mu \leq 1\right\},$$
and the weak $E$-Orlicz-Lorentz space is
$$\Lambda_{\omega,\Phi(E),weak} = \{f\in X_{(0,\infty)}: \|f\|_{\omega,\Phi(E),weak} < \infty\},$$
$$\|f\|_{\omega,\Phi(E),weak}$$
$$= \inf\{\lambda>0: \Phi(E(t,f^*(t)/\lambda))W(t) \leq 1\},$$
$$f^*(t) = \sup\{u: \mu(|f|\geq u) \geq t\}$$
for all $f\in \Lambda_{\omega,\Phi(E)}$.
If $\omega(t)=1$ for $t\in(0,\infty)$, then
$$\Lambda_{\omega,\Phi(E)}(\Omega,\Sigma,\mu) = EL_{\Phi(E)}(\Omega,\Sigma,\mu),$$
$$\Lambda_{\omega,\Phi(E),weak}(\Omega,\Sigma,\mu) = EL_{\Phi(E),weak}(\Omega,\Sigma,\mu).$$
If $\Phi(E(t,u))=u^p$ for $1\leq p<\infty$, we get the Lorentz space
$$\Lambda_{\omega,\Phi(E)}(\Omega,\Sigma,\mu) = EL_{\omega,p}(\Omega,\Sigma,\mu),$$
and the weak Lorentz space
$$\Lambda_{\omega,\Phi(E),weak}(\Omega,\Sigma,\mu) = EL_{\omega,p,weak}(\Omega,\Sigma,\mu).$$
And if $\omega(t)=1$ for $t\in(0,\infty)$, and $\Phi(E(t,u))=u^p$ for $1\leq p<\infty$, then
$$\Lambda_{\omega,\Phi(E)}(\Omega,\Sigma,\mu) = EL_p(\Omega,\Sigma,\mu),$$
$$\Lambda_{\omega,\Phi(E),weak}(\Omega,\Sigma,\mu) = EL_{p,weak}(\Omega,\Sigma,\mu).$$

**Theorem 79.** If $\Phi:\Omega\times[0,\infty)\to\mathbb{R}$ is an increasing $E$-Young function with respect to $E_1,E_2:\Omega\times[0,\infty)\to\Omega\times[0,\infty)$ such that, for $\mu$-a.e. $t\in\Omega$, $E_1(t,x)\leq E_2(t,x)$. Then
$$\Lambda_{\omega,\Phi(E_2)}(\Omega,\Sigma,\mu) \subseteq \Lambda_{\omega,\Phi(E_1)}(\Omega,\Sigma,\mu)$$
and
$$\Lambda_{\omega,\Phi(E_2),weak}(\Omega,\Sigma,\mu) \subseteq \Lambda_{\omega,\Phi(E_1),weak}(\Omega,\Sigma,\mu).$$

**Theorem 80.** If $\Phi_1,\Phi_2:\Omega\times[0,\infty)\to\mathbb{R}$ are $E$-Young functions with respect to $E:\Omega\times[0,\infty)\to\Omega\times[0,\infty)$ such that, for $\mu$-a.e. $t\in\Omega$, $\Phi_1(E(t,x))\leq \Phi_2(E(t,x))$. Then
$$\Lambda_{\omega,\Phi_2(E)}(\Omega,\Sigma,\mu) \subseteq \Lambda_{\omega,\Phi_1(E)}(\Omega,\Sigma,\mu)$$
and
$$\Lambda_{\omega,\Phi_2(E),weak}(\Omega,\Sigma,\mu) \subseteq \Lambda_{\omega,\Phi_1(E),weak}(\Omega,\Sigma,\mu).$$

## VI. CONCLUSION

We have shown that the non $N$-functions, non Young functions, non strong Young functions and non Orlicz functions can be transferred using the $E$-convex theory to $E$-$N$-functions, $E$-Young functions, $E$-strong Young functions and $E$-Orlicz functions respectively. We also have shown that the Orlicz spaces can be generated by non-Young functions but $E$-Young functions with an appropriate map $E$ to extend and generalize studying the Orlicz theory. Moreover, we have considered the inclusion properties of $E$-Orlicz spaces based on effects of the map $E$.